\documentclass[12pt, intlimits]{amsart}
\usepackage[cp1251]{inputenc}
\usepackage[english]{babel}
\usepackage{amsmath, amssymb, amsthm, esint}
\usepackage{geometry}
\usepackage{enumitem}
\usepackage{hyperref}
\usepackage{comment}

\newcommand{\Tr}{\operatorname{Tr}}
\newcommand{\tr}{\operatorname{tr}}

\newcommand{\ind}{\operatorname{ind}}
\newcommand{\Op}{\operatorname{Op}}

\newcommand{\lra}{\longrightarrow}

\newcommand{\CC}{\mathbb{C}}

\newcommand{\RR}{\mathbb{R}}
\renewcommand{\SS}{\mathbb{S}}

\newcommand{\ZZ}{\mathbb{Z}}
\newcommand{\NN}{\mathbb{N}}

\newcommand{\bfa}{\mathbf{a}}
\newcommand{\Mat}[1]{\mathrm{Mat}_{#1}}
\newcommand{\Av}{\operatorname{Av}}
\newcommand{\indt}{\operatorname{ind_t}}
\newcommand{\AS}{\operatorname{AS}}

\newcommand{\I}[1]{\operatorname{I}_{#1}(a)}
\newcommand{\II}{\operatorname{II}(\bfa)}

\newtheorem{theorem}{Theorem}[section]

\newtheorem{proposition}[theorem]{Proposition}

\theoremstyle{definition}
\newtheorem{definition}[theorem]{Definition}
\newtheorem{remark}[theorem]{Remark}

\title[Index of two-dimensional elliptic end-periodic operators]{On the index of two-dimensional elliptic end-periodic operators}
\author{A.Yu. Savin, K.N. Zhuikov}
\email{a.yu.savin@gmail.com, zhuykovcon@gmail.com}

\date{}

\address{RUDN University, Moscow, Russia}

\thanks{MSC classes: 58J20, 58J28, 58J40, 19K56}
\thanks{ArXiv Subjects: Differential Geometry (math.DG); Analysis of PDEs (math.AP); K-Theory and Homology (math.KT)}
\thanks{The work was supported by the Ministry of Science and Higher Education of the Russian Federation within the state assignment (project No. FSSF-2026-0013).}

\begin{document}


\begin{abstract}
This paper investigates the index problem for elliptic end-periodic differential operators on a two-dimensional manifold with a cylindrical end.
We define the topological index as the sum of the contribution of the interior symbol, expressed by an Atiyah--Singer type integral, and the contribution of the limit operator, expressed in terms of the Melrose type eta-invariant and a correction term that depends on the lower order terms of the limit operator.
We prove the homotopy invariance of the topological index, which allows us to reduce the original problem to the study of operators of a simpler form by employing $K$-theoretic methods.
\end{abstract}


\maketitle

\section*{Introduction}

The index problem for elliptic operators was posed by Gelfand~\cite{Gel1} and consists in expressing the analytical index of an elliptic operator in terms of topological invariants.
A solution to this problem for smooth closed manifolds was provided by Atiyah and Singer~\cite{AtSi0}.
This result subsequently opened up an entire field of study at the intersection of differential geometry, topology, and analysis, giving rise to numerous generalizations.

A natural step in the development of this field was the extension of index theory to noncompact manifolds, including manifolds with cylindrical ends.
Such manifolds are defined as noncompact manifolds that are isomorphic to a semi-infinite cylinder outside a compact set.
It should be noted that elliptic theory on manifolds with cylindrical ends is closely related to elliptic theory on manifolds with conical singularities (see, e.g., \cite{NSScS99}).
In the series of papers by Atiyah, Patodi, and Singer~\cite{APS0}, the index problem for Dirac-type operators on such manifolds was investigated.
A key innovation in these works was the introduction of the $\eta$-invariant.
It is a spectral invariant that measures the asymmetry of an operator's spectrum and appears in the index formula as a correction term accounting for the contribution from infinity.
An important generalization of the Atiyah--Patodi--Singer $\eta$-invariant was discovered by Melrose~\cite{Mel6}.
Using the theory of parameter-dependent pseudodifferential operators, Melrose extended the concept of the winding number of a parameter-dependent family determinant to the infinite-dimensional case.
Subsequently, the index problem on manifolds with singularities was studied by many authors; see, e.g.,~\cite{NSScS99,Mel2,AlGe1,Zhu1,AbSaZh2} and the literature cited therein.

In~\cite{FST3,FST6}, Fedosov, Schulze, and Tarkhanov derived an index formula for two-dimensional elliptic pseudodifferential operators.
In this formula, the interior contribution reduces to the standard Atiyah--Singer integral, while the contribution from infinity consists of two terms: one is the $\eta$-invariant defined by the conormal symbol of the operator on the base of a cylinder, and the other is expressed in terms of the principal and subprincipal symbols of the operator.

In recent years, there has been significant development in the index theory for elliptic end-periodic operators (operators with coefficients periodic at infinity).
In the general case, it is natural to construct such a theory on manifolds with periodic ends (end-periodic manifolds), which serve as a natural generalization of manifolds with cylindrical ends.
Manifolds of this type arise, for instance, in geometry and topology in connection with nonstandard smooth structures on $\mathbb{R}^4$~\cite{Tau1} and singular Yamabe metrics on manifolds with singularities~\cite{MPU1}.
In the one-dimensional case, the index problem has been investigated in~\cite{BKS1,Melo1,BoCa08,InRi1}.

In~\cite{MRS2} (for further development of this approach, see~\cite{HaMa1,Tay1}), the index problem for Dirac-type operators on end-periodic manifolds was solved.
In particular, a periodic eta-invariant was defined, which generalizes the Atiyah--Patodi--Singer eta-invariant.
An index formula for end-periodic operators of arbitrary order on the real line was obtained by the authors of the present paper in~\cite{SaZh2}, where Melrose's approach to defining the $\eta$-invariant was used.
The work initiated in the aforementioned paper was continued in~\cite{AfSaZh}, where the periodic $\eta$-invariant was constructed for the case of an infinite cylinder whose base is a smooth closed manifold of arbitrary dimension.

The aim of this paper is to establish the index theorem for elliptic end-periodic operators on a two-dimensional manifold with a cylindrical end.
Let us briefly outline the contents of the work.
Section~\ref{sec:1} introduces the fundamental concepts and statement of the problem.
In Section~\ref{sec:2}, the symbol of the operator under consideration, the $\eta$-invariant, and the topological index are defined, and the main result is stated --- a theorem on the equality of the topological and analytical indices for differential operators.
Section 3 is devoted to the proof of the main result, which is carried out in several steps.
First, the problem is stated in the more general class of pseudodifferential operators.
Second, the homotopy invariance of the topological index is proved.
Finally, it is shown that there exists a homotopy to operators of a simpler form, for which the index formula is already known.
Hence follows the main result.

\section{Statement of the problem}\label{sec:1}

Let $\mathcal{M}$ denote a smooth noncompact two-dimensional manifold diffeomorphic to a semi-infinite cylinder $\mathbb{S}^1 \times (-\infty, 2\pi)$ outside a compact set.
Denote the local coordinates by $x \in \mathbb{S}^1$ and $t \in (-\infty, 2\pi)$, and the dual coordinates by $\xi$ and $p$, respectively.
In what follows, we identify the semi-infinite cylinder $\mathbb{S}^1 \times (-\infty, 2\pi)$ with the corresponding neighborhood $U$ in $\mathcal{M}$.
Let $M$ denote the smooth compact manifold
\[
M = \mathcal{M} \setminus \bigl((-\infty, 0) \times \SS^1\bigr)
\]
with boundary $\partial M = \SS^1 \times \{0\}$.
The noncompact manifold $\mathcal{M}$ is called a \textit{manifold with a cylindrical end}.

Consider the following differential operator of order $\leq m$ on $\mathcal{M}$:
\begin{equation}\label{eq:mainop}
A \colon H^s(\mathcal{M}, E) \lra H^{s-m}(\mathcal{M}, F),
\end{equation}
where $E$ and $F$ are vector bundles over $\mathcal{M}$ with trivializations
\[
E|_U \simeq F|_U \simeq U \times \CC^N,
\]
and $H^s(\cdot)$ is the Sobolev space.
Operator $A$ restricted to $U$ has the form
\begin{equation}\label{eq:DO}
A|_U = \sum_{\alpha + \beta \leq m} a_{\alpha \beta}(x, t)(-i\partial_x)^\alpha(-i\partial_t)^\beta,
\end{equation}
where $\alpha$ and $\beta$ are nonnegative integers.
From here on, the notation $\partial_x = \partial/\partial x$ is used.
We assume that the coefficients $a_{\alpha\beta}(x, t) \in C^\infty\bigl(U, \Mat{N}(\mathbb{C})\bigr)$ are $2\pi$-periodic in $t$ for $t \leq 2\pi$.

\begin{definition}
Operator~\eqref{eq:mainop} is called \textit{elliptic} if
\begin{itemize}
\item[1)]
its \textit{interior symbol}
\[
a = \sum_{\alpha + \beta = m}a_{\alpha\beta}(x,t)\xi^\alpha p^\beta \colon \pi^*E \lra \pi^*F
\]
is an isomorphism outside the zero section in $T^*\mathcal{M}$, where $\pi^*E$ and $\pi^*F$ are the induced bundles corresponding to the projection $\pi \colon T^*\mathcal{M} \to \mathcal{M}$ of the cotangent bundle $T^*\mathcal{M}$ of $\mathcal{M}$;
\item[2)]
the corresponding \textit{limit operator} at infinity is elliptic and invertible:
\begin{equation}\label{eq:limop}
\bfa = \sum_{\alpha + \beta \leq m} a_{\alpha\beta}^\infty(x, t)(-i\partial_x)^\alpha(-i\partial_t)^\beta \colon H^s(\SS^1 \times \RR, \CC^N) \lra H^{s-m}(\SS^1 \times \RR, \CC^N),
\end{equation}
where the coefficients
\[
a_{\alpha\beta}^\infty(x,t) = \lim_{k \to +\infty} a_{\alpha \beta}(x, t-2\pi k)
\]
are $2\pi$-periodic in  $t$.
\end{itemize}
\end{definition}

As is known~\cite{Rab1,Tau1}, the ellipticity of operator~\eqref{eq:mainop} implies its fredholmness.
Our goal is to derive an index formula for elliptic operator~\eqref{eq:mainop}.

\begin{remark}
The invertibility of the limit operator~\eqref{eq:limop} is equivalent to the invertibility of the corresponding family of operators $\bfa(\theta)$ on the torus $\SS^1 \times \SS^1$, obtained via the Fourier--Laplace transform, see~\cite[Theorem~2]{Rab1} and~\cite[Lemma~4.3]{Tau1}.
Specifically, conjugating the periodic operator $\bfa$ with the Fourier--Laplace transform $\mathcal{F}_{\theta} \colon L^2(\RR_t) \to L^2(\SS^1_t \times \SS^1_\theta)$ yields the family
\begin{equation}\label{eq:atheta}
\begin{gathered}
\bfa(\theta) = \mathcal{F}_\theta \bfa \mathcal{F}_\theta^{-1} \colon H^s(\SS^1_x \times \SS^1_t, \CC^N) \otimes L^2(\SS^1_\theta) \lra H^{s-m}(\SS^1_x \times \SS^1_t, \CC^N) \otimes L^2(\SS^1_\theta), \\
\bfa(\theta) = \sum_{\alpha + \beta \leq m} a_{\alpha\beta}^\infty(x, t)(-i\partial_x)^\alpha\Bigl(-i\partial_t-\frac{\theta}{2\pi}\Bigr)^\beta.
\end{gathered}
\end{equation}
We call family~\eqref{eq:atheta} the \textit{Gelfand transform} of $\bfa$.
\end{remark}

\section{Index formula}\label{sec:2}
Consider an elliptic differential operator of the form~\eqref{eq:mainop}.
\subsection{Contribution of the interior symbol}
Denote the interior symbol of operator~\eqref{eq:mainop} by
\[
a \in C^\infty\bigl(S^*\mathcal{M}, \mathrm{Hom}(\pi^*E, \pi^*F)\bigr),
\]
where $S^*\mathcal{M} = \{(x, t, \xi, p) \in T^*\mathcal{M} \mid \xi^2 + p^2 = 1\}$ is the cosphere bundle of $\mathcal{M}$ with respect to a Riemannian metric $g$ on $\mathcal{M}$, which takes the form $g = dx^2 + dt^2$ outside the submanifold $M$.

Let $\nabla_E$ and $\nabla_F$ denote connections on the vector bundles $E$ and $F$, and let $\theta_E = \nabla_E^2$ and $\theta_F = \nabla_F^2$ be the curvatures of these connections.
We assume that
\[
\nabla_E|_U = \nabla_F|_U = d \text{ is the trivial connection.}
\]
Recall that for any homomorphism $b \colon E \to F$, a covariant differential $\nabla b \colon E \to F \otimes \Lambda^1(\mathcal{M})$ is defined by the equality
\[
(\nabla b) u = \nabla_F(bu) - b\nabla_Eu,
\]
where $u \in C^\infty(\mathcal{M},E)$ is an arbitrary section.

\begin{definition}
For an elliptic interior symbol $a$, the \textit{Atiyah--Singer form} is defined by the expression (see~\cite[Chap.~1, Sec.~4]{Fds2})
\[
\AS(a) = \frac{1}{24\pi^2}\tr(a^{-1}\nabla a)^3 - \frac{1}{8\pi^2}\tr\bigl(a^{-1}\nabla a(\theta_E + a^{-1}\theta _F a)\bigr) \in \Omega^3(S^*\mathcal{M}),
\]
where $\Omega^3(S^*\mathcal{M})$ is the space of differential forms of degree $3$ on $S^*\mathcal{M}$.

The \textit{contribution of the interior symbol} to the index formula is defined as
\begin{equation}\label{eq:I}
\I{} = \int_{S^*M}\AS(a) + \int_{S^*([-2\pi, 0] \times \SS^1)}\hspace{-.5cm}\frac{t}{2\pi}\AS(a).
\end{equation}
\end{definition}

\begin{remark}\label{rem:AScyl}
At the cylindrical end (for $t \leq 2\pi$), $\nabla=d$ and $\theta_E = \theta_F = 0$.
Therefore,
\begin{equation}\label{eq:AScyl}
\AS(a) = \frac{1}{24\pi^2}\tr(a^{-1}da)^3.
\end{equation}
\end{remark}

\begin{remark}\label{rem:t0}
We claim that for all $t_0 < 2\pi$, the following equality holds:
\begin{equation}\label{eq:It0}
\I{} = \int_{S^*(\mathcal{M}\setminus(-\infty, t_0) \times \SS^1)}\hspace{-1cm}\AS(a) + \int_{S^*([t_0-2\pi, t_0] \times \SS^1)}\hspace{-.5cm}\frac{t}{2\pi}\AS(a).
\end{equation}
That is, the choice of the cut location $t_0$ on the cylinder does not affect the value of the functional on the right-hand side in~\eqref{eq:It0}.
Indeed, denote the right-hand side in~\eqref{eq:It0} by $I_{t_0}(a)$.
Then, for any $t_1 < t_2 < 2\pi$, a direct calculation yields the equalities
\begin{align*}
\I{t_0}\Big|^{t_1}_{t_2} &= \int_{S^*([t_1, t_2] \times \SS^1)}\hspace{-0.5cm}\AS(a) + \int_{S^*([t_1 - 2\pi, t_2 - 2\pi] \times \SS^1)}\hspace{-1cm}\frac{t}{2\pi}\AS(a) - \int_{S^*([t_1, t_2] \times \SS^1)}\hspace{-.5cm}\frac{t}{2\pi}\AS(a) = \\
&= \int_{S^*([t_1, t_2] \times \SS^1)}\hspace{-0.5cm}\left(1 - \frac{t}{2\pi} \right)\AS(a) + \int_{S^*([t_1, t_2] \times \SS^1)}\hspace{-0.5cm} \left(-1 + \frac{t}{2\pi}\right)\AS(a) = 0.
\end{align*}
\end{remark}

\begin{remark}
Note that the integral of the Atiyah--Singer form $\AS(a)$ over a period is equal to zero.
Indeed, it is equal to the index of the operator $\bfa(\theta)$ on the torus, and this operator is invertible by the ellipticity condition.
\end{remark}

\subsection{Eta-invariant. Contribution of the limit operator at infinity}

By the ellipticity of operator $A$, the limit operator $\mathbf{a}$ is invertible.
Consider the operator
\[
\bfa^{-1}[t, \bfa] \colon H^s(\SS^1 \times \RR, \CC^N) \lra H^{s}(\SS^1 \times \RR, \CC^N),
\]
where $[\cdot, \cdot]$ denotes a commutator.
This operator is invariant under shifts by $2\pi$.
Its averaging over the period
\begin{equation}\label{eq:Avs}
\Av\bigl(\bfa^{-1}[t, \bfa]\bigr)
\end{equation}
is a translation-invariant operator on the cylinder $S^1 \times \mathbb{R}$.
Here the averaging of a periodic operator is defined by the formula
\[
B \longmapsto \Av(B) = \frac{1}{2\pi}\int_0^{2\pi}T_\varphi B T_{-\varphi}d\varphi, \quad \text{where}\quad T_\varphi u(x,t) = u(x, t + \varphi).
\]
One can show (see~\cite[Sec.~1]{AfSaZh}) that operator~\eqref{eq:Avs} is a translation-invariant pseudodifferential operator obtained from a parameter-dependent pseudodifferential operator using the Fourier transform $\mathcal{F}_{t \to p}$:
\[
A(p) \longmapsto A(-i\partial_t) = \mathcal{F}^{-1}_{p \to t}A(p)\mathcal{F}_{t \to p}.
\]

Let $\Tr$ be Melrose's regularized trace~\cite{Mel6} on translation-invariant pseudodifferential operators.
We define a functional $\Tr_1$ by the expression
\[
\Tr_1 \bfa = \Tr \circ \Av [it, \bfa] = \frac{1}{4\pi^2} \int_{T^*\SS^1 \times \SS^1}\tr a_1(x, t, \xi, p)\Big|^{p=1}_{p=-1}dxd\xi dt,
\]
where $a_1$ is the component of order $-1$ of the complete symbol of $\bfa$.

\begin{definition}
For an invertible limit operator $\bfa$ of order $m$, the \textit{$\eta$-invariant} is defined as the number
\begin{equation}\label{eq:overeta}
\eta(\bfa) = \frac{1}{2\pi}\Tr\circ\Av\bigl(\bfa^{-1}[t, \bfa]\bigr) \in \CC.
\end{equation}
The \textit{contribution of the limit operator} to the index formula is defined as
\begin{equation}\label{eq:II}
\II = -\eta(\bfa) + \frac{1}{4\pi}\Tr_1\bigl(a^{-1}(a_\xi a^{-1}a_x - 2ia_1)\bigr), \end{equation}
where $a$ and $a_1$ are the components of the complete symbol of $\mathbf{a}$ of orders $m$ and $m-1$, respectively, and the subscripts denote partial derivatives with respect to the corresponding variables.
\end{definition}

\begin{definition}
The \textit{topological index} of an elliptic operator $A$ is defined as the sum of the contributions from its interior symbol~\eqref{eq:I} and the limit operator~\eqref{eq:II}:
\begin{equation}\label{eq:indt}
\indt A = \I{} + \II.
\end{equation}
\end{definition}

The following theorem is the main result of this work.

\begin{theorem}\label{th:main}
Let $A \colon H^s(\mathcal{M}, E) \to H^{s-m}(\mathcal{M}, F)$ be an elliptic end-periodic differential operator of order $m$.
Then
\begin{itemize}
\item[1)]
its topological index~\eqref{eq:indt} is invariant under homotopies of $A$ within the class of elliptic operators;
\item[2)]
the following index formula holds:
\begin{equation}\label{eq:find}
\ind A = \indt A,
\end{equation}
where the left-hand side of the equality contains the Fredholm index.
\end{itemize}
\end{theorem}
The proof of Theorem~\ref{th:main} is given in the next section.
Here we present a brief sketch of the proof.

First, for a family of elliptic interior symbols $a = a(\varepsilon)$, we compute the derivative $\partial_\varepsilon \I{}$ with respect to the parameter $\varepsilon$.
Next, for a family of end-periodic operators $\bfa = \bfa(\varepsilon)$, we compute the derivative $\partial_\varepsilon \II$.
We then show that, for an elliptic family $A(\varepsilon)$, these derivatives cancel each other out.
This establishes the homotopy invariance stated in part 1 of the theorem.
To prove part 2 of the theorem, we provide a homotopy classification of elliptic operators.
This enables us to reduce the problem to operators of a simpler form, for which it is straightforward to verify that the left- and right-hand sides in~\eqref{eq:find} are equal.

\section{Proof of the main result}\label{sec:3}
We will obtain Theorem~\ref{th:main} as a special case of the index theorem for the corresponding pseudodifferential operators ($\psi$DOs).
Let us define these $\psi$DOs and state the theorem.

\subsection{Index theorem for $\psi$DOs}
The \textit{space of interior principal symbols} is the space
\[
\Sigma_{\mathrm{int}}=\Bigl\{a \in C^\infty\bigl(S^*\mathcal{M}, \mathrm{Hom}(\pi^*E, \pi^*F)\bigr) \mid a(x, t, \xi, p) = a(x, t - 2\pi k, \xi, p)\, \forall t \leq 2\pi, k \in \NN \Bigr\}.
\]
The \textit{algebra of periodic $\psi$DOs of order zero on the cylinder $\SS^1 \times \mathbb{R}$} is the algebra (see~\cite[Sec.~1]{AfSaZh})
\[
\Psi^0_{\mathrm{per}}(\SS^1 \times \RR) = \Bigl\{\bfa = \sum_{k \in \ZZ}e^{ikt}A_k(x, -i\partial_x, -i\partial_t) \colon H^s(\SS^1 \times \RR, \CC^N) \lra H^{s}(\SS^1 \times \RR, \CC^N)\Bigr\},
\]
where $A_k = \mathcal{F}^{-1}_{p \to t} A_k(p)\mathcal{F}_{t \to p}$ and $A_k(p)$ is a parameter-dependent $\psi$DO of order $\leq 0$ on $\SS^1$ with parameter $p \in \RR$.
Here the coefficients $A_k$ are assumed to decay rapidly as $k \to \infty$.

The \textit{complete symbol} of an operator $\bfa \in \Psi^0_{\mathrm{per}}(\SS^1 \times \RR)$ is the formal series
\begin{equation}\label{eq:compsym}
\sigma(\bfa) = \sum_{k \in \ZZ} e^{ikt}a_{k0} + \sum_{k \in \ZZ} e^{ikt}a_{k1} + \dotso,
\end{equation}
where $a_{kj}$ is the $j$-th component of the complete symbol of operator $A_k$.
The leading component
\[
\sigma_0(\bfa) = \sum_{k \in \ZZ}e^{ikt}a_{k0} \in \Mat{N}\bigl(C^\infty(S^*(\SS^1 \times \SS^1))\bigr)
\]
of the complete symbol~\eqref{eq:compsym} is called the \textit{principal symbol} of operator $\bfa$.

\begin{definition}
The \textit{space of principal symbols} is defined by the formula
\[
\Sigma = \Big\{(a, \bfa) \in \Sigma_{\mathrm{int}} \oplus \Psi^0_{\mathrm{per}}(\SS^1 \times \RR) \;\Big|\; a|_{\SS^1 \times [0, 2\pi]} = \sigma_0(\bfa)\Big\}
\]
and consists of interior symbols $a$ and limit operators $\mathbf{a}$ subject to compatibility condition.
For brevity, we refer to these principal symbols simply as symbols.
\end{definition}

We define the quantization mapping
\begin{equation} \label{eq:pdo}
\Op \colon \Sigma \lra \mathcal{B}\bigl(L^2(\mathcal{M}, E), L^2(\mathcal{M}, F)\bigr), \quad \Op (a, \bfa) = \psi_0\Op(a)\chi_0 + \psi_1\bfa\chi_1,
\end{equation}
where $\Op(a)$ is a $\psi$DO with principal symbol $a$, $\chi_0 + \chi_1 = 1$ is a partition of unity on $\mathcal{M}$, the function $\chi_0$ is equal to $1$ on $M$ and $0$ on $\SS^1 \times (-\infty, -\pi)$, and the support of the function $\psi_j$ lies within a small neighborhood of the support of the function $\chi_j$, with $\psi_j \chi_j = \chi_j$.

A symbol $(a, \mathbf{a}) \in \Sigma$ is called \textit{elliptic} if there exists the inverse symbol $(a^{-1}, \mathbf{a}^{-1}) \in \Sigma$ corresponding to the bundles $F$ and $E$.
As in the case of differential operators, the ellipticity of $(a, \mathbf{a})$ implies the fredholmness of the operator $\Op(a, \mathbf{a})$.
In this case, $\Op(a^{-1}, \mathbf{a}^{-1})$ is an almost inverse of $\Op(a, \mathbf{a})$:
\[
\Op(a^{-1}, \bfa^{-1}) \Op(a, \bfa) = 1 + K_1, \quad \Op(a, \bfa) \Op(a^{-1}, \bfa^{-1}) = 1 + K_2,
\]
where $K_j$ are compact operators.

Let us state an analogue of Theorem~\ref{th:main} for elliptic $\psi$DOs.

\begin{theorem}\label{th:main'}
Let $(a, \bfa) \in \Sigma$ be an elliptic symbol. Then
\begin{itemize}
\item[1)]
the expression $\indt(a, \bfa) = \I{} + \II$ {\rm (}cf.~\eqref{eq:indt}{\rm )} is invariant under homotopies within the class of elliptic symbols;
\item[2)]
the following index formula holds true:
\[
\ind \Op(a, \bfa) = \indt(a, \bfa).
\]
\end{itemize}
\end{theorem}

\subsection{Homotopy invariance}
Let us prove part 1 of Theorem~\ref{th:main'}.
The proof is divided into three steps.
First, we derive a formula for the variation of the contribution of the interior symbol.
Next, we obtain a formula for the variation of the contribution of the limit operators resulting from changes in the lower order terms.
Finally, we arrive at a general formula for the variation of the contribution of the limit operators.

\begin{proposition}\label{prop:dsI}
Let $a = a(s) \in \Sigma_{\mathrm{int}}$, $s \in [0,1]$, be a smooth homotopy of elliptic interior symbols.
Then
\begin{equation}\label{eq:dsI}
\partial_s \I{} = \frac{1}{4\pi}\Tr_1\bigl(a^{-1}a_s[a^{-1}a_x, a^{-1}a_\xi]\bigr).
\end{equation}
\end{proposition}
\begin{proof}
1. First, we prove the equality
\begin{equation}\label{eq:dsI1}
\partial_s \I{} = -\frac{1}{16\pi^3}\int_{S^*([0,2\pi] \times \SS^1)} dt \tr\bigl(a^{-1}a_s(a^{-1} da)^2\bigr).
\end{equation}
To this end, we calculate the derivative $\partial_s \AS(a)$:
\begin{equation}\label{eq:dsas}
\partial_s\AS(a) = -\frac{1}{8 \pi^2} \partial_s \tr (a^{-1}\nabla a\theta_E)  -\frac{1}{8 \pi^2} \partial_s \tr (\nabla aa^{-1}\theta_F) + \frac{1}{24 \pi^2} \partial_s \tr (a^{-1}\nabla a)^3.
\end{equation}
For the first two terms on the right-hand side in~\eqref{eq:dsas}, a straightforward calculation yields
\begin{equation}\label{eq:dsas1}
\begin{aligned}
\partial_s \tr (a^{-1}\nabla a\theta_E) &+ \partial_s \tr (\nabla aa^{-1}\theta_F)  \\
=& -\tr(a^{-1}a_sa^{-1}\nabla a \theta_E - a^{-1}\nabla a_s\theta_E) + \tr(\nabla a_s a^{-1}\theta_F - \nabla a a^{-1}a_sa^{-1}\theta_F)  \\
=& -\tr\bigl(a^{-1}a_sa^{-1}\nabla a \theta_E - \nabla(a^{-1}a_s\theta_E) - a^{-1}\nabla a a^{-1} a_s \theta_E\bigr)  \\
&+\tr\bigl(\nabla(a_s a^{-1} \theta_F) + a_s a^{-1} \nabla a a^{-1} \theta_F - a^{-1}a_sa^{-1}\theta_F \nabla a\bigr).
\end{aligned}
\end{equation}
For the last term in~\eqref{eq:dsas}, we obtain (cf.~\cite[Lemma~1.17]{Zha3})
\begin{equation}\label{eq:dsas2}
\begin{aligned}
&\partial_s \tr (a^{-1}\nabla a)^3 = 3 \tr \bigl(\partial_s(a^{-1}\nabla a)(a^{-1} \nabla a)^2\bigr)  \\
=& 3 \tr\bigl([a^{-1} \nabla a, a^{-1} a_s](a^{-1} \nabla a)^2 \bigr) + 3 \tr \bigl(\nabla(a^{-1} a_s)(a^{-1}\nabla a)^2\bigr)  \\
=& 3 \tr [a^{-1}\nabla a, a^{-1} a_s (a^{-1}\nabla a)^2] + 3 \tr \Bigl(\nabla \bigl(a^{-1} a_s (a^{-1}\nabla a)^2\bigr)\Bigr) - 3 \tr \bigl(a^{-1} a_s \nabla (a^{-1}\nabla a)^2\bigr) \\
=& 3 \tr \Bigl(\nabla \bigl(a^{-1} a_s (a^{-1}\nabla a)^2\bigr)\Bigr) - 3 \tr\bigl(a^{-1} a_s (a^{-1}\theta_F \nabla a - \theta_E a^{-1} \nabla a - a^{-1} \nabla a a^{-1} \theta_F a + a^{-1} \nabla a \theta_E)\bigr).
\end{aligned}
\end{equation}
Here, trace of the commutator in the last equality vanishes by the cyclic property of the trace.

Substituting~\eqref{eq:dsas1} and~\eqref{eq:dsas2} into~\eqref{eq:dsas}, we obtain
\begin{equation}\label{eq:dsasfin}
\begin{aligned}
\partial_s \AS(a) &= \frac{1}{8\pi ^2} \tr \nabla \bigl(a^{-1} a_s (a^{-1}\nabla a)^2\bigr) -\frac{1}{8\pi ^2} \tr \bigl(\nabla (a^{-1}a_s \theta_E) + \nabla (a_sa^{-1}\theta_F)\bigr)  \\
&= \frac{1}{8\pi ^2} d \tr \bigl (a^{-1} a_s (a^{-1}\nabla a)^2\bigr) -\frac{1}{8\pi ^2} d \tr \bigl ( a^{-1}a_s \theta_E + a_sa^{-1}\theta_F\bigr).
\end{aligned}
\end{equation}
Given Remark~\ref{rem:AScyl}, we integrate~\eqref{eq:dsasfin} and apply Stokes' theorem:
\begin{equation}\label{eq:intdsas}
\begin{aligned}
\int_{S^*M}\partial_s \AS(a) &= \frac{1}{8\pi^2}\int_{S^*M} d \tr \bigl (a^{-1} a_s (a^{-1}\nabla a)^2\bigr) -\frac{1}{8\pi^2}\int_{S^*M} d \tr \bigl ( a^{-1}a_s \theta_E + a_sa^{-1}\theta_F\bigr)  \\
&= \frac{1}{8\pi^2}\int_{\partial S^*M} \tr \bigl (a^{-1} a_s (a^{-1}d a)^2\bigr).
\end{aligned}
\end{equation}

Now, we calculate the derivative of the second term in~\eqref{eq:I} using~\eqref{eq:dsasfin}:
\begin{equation}\label{eq:dstas}
\begin{aligned}
\partial_s \hspace{-0.5cm} \int_{S^*([0, 2\pi] \times \SS^1)}\hspace{-0.5cm}\frac{t}{2\pi}\AS(a) =& \frac{1}{8\pi^2} \hspace{-0.5cm} \int_{S^*([0, 2\pi] \times \SS^1)}\hspace{-0.5cm}\frac{t}{2\pi}d\tr\bigl(a^{-1}a_s(a^{-1}da)^2\bigr)  \\
=& \frac{1}{8\pi^2} \hspace{-0.5cm} \int_{S^*([0,2\pi] \times \SS^1) \cap \{t = 2\pi\}}\hspace{-.5cm}\tr\bigl(a^{-1}a_s(a^{-1}da)^2\bigr) - \frac{1}{16\pi^3} \hspace{-0.5cm} \int_{S^*([0, 2\pi] \times \SS^1)}\hspace{-.5cm}dt\tr\bigl(a^{-1}a_s(a^{-1}da)^2\bigr)  \\
=& -\frac{1}{8\pi^2}\int_{\partial S^*M}\tr\bigl(a^{-1}a_s(a^{-1}da)^2\bigr) - \frac{1}{16\pi^3} \hspace{-0.5cm} \int_{S^*([0, 2\pi] \times \SS^1)}\hspace{-.5cm}dt\tr\bigl(a^{-1}a_s(a^{-1}da)^2\bigr).
\end{aligned}
\end{equation}
Here we used Remark~\ref{rem:AScyl} and the integration-by-parts formula in the second equality, and the periodicity of the interior symbol $a$ in the third one.

Finally, adding~\eqref{eq:intdsas} and~\eqref{eq:dstas}, we obtain the desired equality~\eqref{eq:dsI1}.

2. Now the desired equality~\eqref{eq:dsI} will follow from~\eqref{eq:dsI1} if we establish the equality
\[
\int_{S^*([0, 2\pi] \times\SS^1)}\hspace{-5mm}dt \tr\bigl(a^{-1}a_s(a^{-1} da)^2\bigr) = -4\pi^2\Tr_1\bigl(a^{-1}a_s[a^{-1}a_x, a^{-1}a_\xi]\bigr).
\]
Since $S^*([0, 2\pi] \times \SS^1_x) = [0, 2\pi] \times \SS^1_x \times \SS^1_{\xi, p}$ (where $\xi^2 + p^2 = 1$) and the function $a(x, t, \xi, p)$ is homogeneous of degree $0$ in the variables $\xi$ and $p$, the following relations hold:
\[
a_p = -\frac{\xi}{p}a_\xi, \quad dp = -\frac{\xi}{p}d\xi.
\]
These relations imply
\begin{equation}\label{eq:ada2}
\begin{aligned}
\int_{S^*([0, 2\pi] \times \SS^1)}\hspace{-5mm} dt  \tr\bigl(a^{-1}a_s(a^{-1}da)^2\bigr) =& -\hspace{-9mm}\int_{[0, 2\pi] \times \SS^1 \times (-1, 1)}\hspace{-11mm} p^{-2}\tr\bigl(a^{-1}a_s[a^{-1}a_x, a^{-1}a_\xi]\bigr)(x, t, \xi, \alpha |p|)\Big|^{\alpha = 1}_{\alpha = -1} dt dx d\xi  \\
=& -\hspace{-9mm}\int_{[0, 2\pi] \times \SS^1 \times (-1, 1)}\hspace{-11mm} p^{-3}\tr\bigl(a^{-1}a_s[a^{-1}a_x, a^{-1}a_\xi]\bigr)\Bigl(x, t, \frac{\xi}{|p|}, \alpha\Bigr)\Big|^{\alpha = 1}_{\alpha = -1} dt dx d\xi  \\
=& -\hspace{-5mm}\int_{T^*\SS^1 \times [0, 2\pi]}\hspace{-5mm} \tr\bigl(a^{-1}a_s[a^{-1}a_x, a^{-1}a_\xi]\bigr)(x, t, \tau, \alpha)\Big|^{\alpha = 1}_{\alpha = -1} dt dx d\tau = \\
=& -4\pi^2\Tr_1\bigl(a^{-1}a_s[a^{-1}a_x, a^{-1}a_\xi]\bigr).
\end{aligned}
\end{equation}
In the second equality, we factored out $p^{-1}$ using the homogeneity of degree $-1$.
In the third equality, we made the substitution $\tau = \xi |p|^{-1}$ (with $d\tau = |p|^{-3}d\xi$).
In the final step, we relabeled $\tau$ and $\alpha$ as $\xi$ and $p$, respectively.
Finally, by substituting this expression into~\eqref{eq:dsI1}, we arrive at the desired equality~\eqref{eq:dsI}.

Proposition~\ref{prop:dsI} is proved.
\end{proof}

\begin{proposition}\label{prop:homprincipal}
Let $\mathbf{a} \in \Psi^0_{\mathrm{per}}(\mathbb{S}^1 \times \mathbb{R})$ be an elliptic invertible $\psi$DO with complete symbol $\sigma(\mathbf{a}) = a + a_1 + \dotsb$, where $a$ is the principal symbol, $a_1$ is the component of order $-1$, and the ellipsis denotes terms of lower order.
Consider a smooth homotopy of elliptic invertible operators $\mathbf{b} = \mathbf{b}(s) \in \Psi^0_{\mathrm{per}}(\mathbb{S}^1 \times \mathbb{R})$, $s \in [0, 1]$, with complete symbol $b = b(s) = a + s(a_1 + \dotsb)$.
Then the following equalities hold:
\begin{gather}\label{eq:detas}
\partial_s \eta(\mathbf{b}) = \frac{1}{2\pi i}\Tr_1 (a^{-1}a_1),\\ \label{eq:detas2}
\partial_s \operatorname{II}(\mathbf{b}) = 0.
\end{gather}
\end{proposition}

\begin{proof}
1.
Let us prove formula~\eqref{eq:detas}.
In~\cite[Theorem 4]{AfSaZh}, the following formula for the derivative of the $\eta$-invariant with respect to the parameter was obtained:
\begin{equation}\label{eq:afsazh}
\partial_s \eta(\mathbf{b}) = \frac{1}{2\pi i}\Tr_1\Bigl(\sigma(\mathbf{b}^{-1}\partial_s\mathbf{b}) - i(a^{-1}a_s)_ta^{-1}a_p \Bigr).
\end{equation}
Hence, for the current situation, we obtain
\begin{equation*}
2\pi i\partial_s\eta(\mathbf{b}) = \Tr_1 (b^{*-1}*b_s).
\end{equation*}
Here we use the composition formula for the complete symbols:
\[
\sigma(\mathbf{b}^{-1}\partial_s\mathbf{b}) = \sigma(\mathbf{b}^{-1}) * \sigma(\partial_s\mathbf{b}) = \sigma(\mathbf{b})^{*-1} * \sigma(\partial_s\mathbf{b}).
\]
By $*$ we denote the Kohn--Nirenberg product:
\[
a*b \sim ab + \sum_{\alpha + \beta \geq 0}\frac{1}{\alpha!\beta!}\bigl[\partial_\xi^\alpha\partial_p^\beta a\bigr] \bigl[(-i\partial_x)^\alpha(-i\partial_t)^\beta b\bigr] = ab -i(a_\xi b_x + a_pb_t) + \dotsi.
\]
The second term in~\eqref{eq:afsazh} vanishes since the principal symbol $a$ does not depend on the parameter $s$ by assumption.

Furthermore, since $b_s = a_1 + \dotsb$, the composition $b^{-1} * b_s$ is of order $-1$, i.e., only its leading component $a^{-1}a_1$ contributes to $\Tr_1$.
Hence, the desired equality~\eqref{eq:detas} follows.

2. Now we prove formula~\eqref{eq:detas2}.
From~\eqref{eq:II} and~\eqref{eq:detas}, we have
\[
\partial_s \operatorname{II}(\mathbf{b}) = -\frac{1}{2\pi i}\Tr_1(a^{-1}a_1) + \frac{1}{2\pi i}\partial_s\Tr_1(a^{-1}sa_1) = 0.
\]
Proposition~\ref{prop:homprincipal} is proved.
\end{proof}

Proposition~\ref{prop:homprincipal} implies that it suffices to consider periodic $\psi$DOs whose complete symbols coincide with their principal symbol.

\begin{proposition}\label{prop:deta}
Let $\bfa = \bfa(s) \in \Psi^0_{\mathrm{per}}(\SS^1 \times \RR)$, $s \in [0,1]$, be a smooth homotopy of invertible elliptic periodic $\psi$DOs of order zero on the cylinder $\SS^1 \times \RR$ with the complete symbol $\sigma(\bfa) = \sigma_0(\bfa) = a(s)$ consisting solely of the principal symbol.
Then
\begin{align}\label{eq:deta}
\partial_s \eta(\bfa) =& -\frac{1}{2\pi}\Tr_1\bigl(a^{-1}a_\xi a^{-1}(a_xa^{-1}a_s - a_{sx})\bigr),\\ \label{eq:dsII}
\partial_s \II =& -\frac{1}{4\pi}\Tr_1\bigl(a^{-1}a_s[a^{-1}a_x, a^{-1}a_\xi]\bigr).
\end{align}
\end{proposition}

\begin{proof}
1. Let us prove formula~\eqref{eq:deta}.
By~\eqref{eq:afsazh}, we obtain
\begin{equation}\label{eq:diffPhi}
2\pi i\partial_s\eta(A) = \Tr_1 \bigl(a^{*-1}*a_s - i(a^{-1}a_s)_ta^{-1}a_p\bigr).
\end{equation}
The first two components of the inverse symbol $a^{*-1} = b_0 + b_1 + \dotsi$ are equal to
\[
b_0 = a^{-1} \quad \text{and} \quad
b_1 = ia^{-1}\bigl(a_\xi(a^{-1})_x + a_p(a^{-1})_t\bigr) = -ia^{-1}(a_\xi a^{-1}a_x + a_p a^{-1}a_t)a^{-1}.
\]
Since only components of order $-1$ contribute to $\Tr_1$, the contribution of the first term $a^{*-1} * a_s$ in~\eqref{eq:diffPhi} is given by
\begin{equation}\label{eq:first}
(a^{-1} * a_s)_1 = ia^{-1}(a_\xi a^{-1} a_{sx} + a_pa^{-1}a_{st}) \quad\text{and} \quad b_1a_s = -ia^{-1}(a_\xi a^{-1}a_x + a_p a^{-1}a_t)a^{-1}a_s.
\end{equation}
The contribution of the second term in~\eqref{eq:diffPhi} is given by the entire expression:
\begin{equation}\label{eq:second}
-i(a^{-1}a_s)_ta^{-1}a_p = -ia^{-1}(a_{st} - a_ta^{-1}a_s)a^{-1}a_p.
\end{equation}
Substituting~\eqref{eq:first} and~\eqref{eq:second} into~\eqref{eq:diffPhi}, we obtain
\begin{equation}\label{eq:detasymb}
\begin{aligned}
2\pi \partial_s\eta(A) &= 
\Tr_1\bigl( a^{-1}a_\xi a^{-1} ( a_{sx} - a_xa^{-1}a_s )\bigr) + \Tr_1[ a^{-1}a_p, a^{-1}a_{st}] + \Tr_1[ a^{-1} a_ta^{-1}a_s, a^{-1}a_p ]  \\
& = \Tr_1\bigl( a^{-1}a_\xi a^{-1} ( a_{sx} - a_xa^{-1}a_s )\bigr).
\end{aligned}
\end{equation}
Here, the last two terms vanish by the equality (see~\cite[Corollary~3, Theorem~3]{AfSaZh})
\begin{equation}\label{eq:Tr2}
\Tr_1[a, b] = \Tr_1(ia_t * b)_p.
\end{equation}
Namely, since the commutators are of order $\leq -2$, they do not contribute to $\Tr_1$.

2. Let us prove formula~\eqref{eq:dsII}.
Substituting~\eqref{eq:detasymb} into~\eqref{eq:II}, we obtain
\begin{equation}\label{eq:dsII2}
\begin{aligned}
\partial_s \II =& -\partial_s \eta(\bfa) + \frac{1}{4\pi} \partial_s \Tr_1(a^{-1} a_\xi a^{-1} a_x ) = \frac{1}{2\pi}\Tr_1(a^{-1}a_\xi a^{-1}a_xa^{-1}a_s - a^{-1}a_\xi a^{-1}a_{sx})  \\
&- \frac{1}{4\pi}\Tr_1(a^{-1} a_s a^{-1} a_\xi a^{-1} a_x - a^{-1} a_{\xi s} a^{-1} a_x + a^{-1} a_\xi a^{-1} a_s a^{-1} a_x -a^{-1} a_\xi a^{-1} a_{xs})  \\
=& -\frac{1}{4\pi}\Tr_1\bigl(a^{-1}a_s[a^{-1}a_x, a^{-1}a_\xi]\bigr) + \frac{1}{4\pi}\Tr_1(a^{-1} a_x a^{-1} a_{\xi s} - a^{-1} a_\xi a^{-1} a_{x s}).
\end{aligned}
\end{equation}
We transform the last term in~\eqref{eq:dsII2}:
\begin{equation}\label{eq:dsII3}
\begin{aligned}
&\hspace{1em}\Tr_1(a^{-1} a_x a^{-1} a_{\xi s} - a^{-1} a_\xi a^{-1} a_{x s})  \\
&= - \Tr_1(i a_s * a^{-1} + a^{-1}a_ta^{-1}a_{ps}) + \Tr_1(i a^{-1} * a_s + a^{-1}a_{p}a^{-1}a_{ts})  \\
&= i\Tr_1[a^{-1},a_s] +  \Tr_1\bigl((a^{-1})_ta_{ps} - (a^{-1})_pa_{ts}\bigr)  \\
&= -\Tr_1\bigl(((a^{-1})_t a_s)_p\bigr) + \Tr_1\bigl((a^{-1})_ta_{ps} - (a^{-1})_pa_{ts}\bigr)  \\
&= -\Tr_1\bigl((a^{-1})_{pt}a_s\bigr) - \Tr_1\bigl((a^{-1})_{p}a_{ts}\bigr) = 0.
\end{aligned}
\end{equation}
Here, in the third equality, we used formula~\eqref{eq:Tr2}, and in the last equality, we performed integration by parts with respect to $t$.
Substituting~\eqref{eq:dsII3} into~\eqref{eq:dsII2}, we obtain the desired equality~\eqref{eq:dsII}.

Proposition~\ref{prop:deta} is proved.
\end{proof}

Finally, the desired homotopy invariance of the topological index follows from Propositions~\ref{prop:dsI},~\ref{prop:homprincipal}, and~\ref{prop:deta}:
\[
\partial_s \indt \bfa(s) = 0.
\]
This completes the proof of part 1 of Theorem~\ref{th:main'}.

\subsection{Homotopy classification}
Let us proceed to the proof of part 2 of Theorem~\ref{th:main'}.
By part 1 of Theorem~\ref{th:main'}, the topological index of an elliptic operator is invariant under homotopies.
We will construct a homotopy from the original operator to an operator of a simpler form, for which the index formula is already known.

\begin{proposition}
\label{prop1}
Let $A$ be an elliptic operator of order zero of the form~\eqref{eq:pdo} with the limit operator $\mathbf{a}$. Then
\begin{enumerate}
\item[{\rm1)}]   there exist integers $N'\ge 1$ and $k_1,k_2,k_3$ and the homotopy 
\begin{equation}
\label{eq:homotopy2}
\mathbf{a}\oplus 1_{N'}\sim e^{i(k_1x+k_2t)}\oplus a (x,-i\partial_x,-i\partial_t)\oplus 1_{N'+N-2} 
\end{equation}
of periodic elliptic invertible operators on the cylinder $\mathbb{S}^1 \times \mathbb{R}$, where $a(x, -i\partial_x, -i\partial_t)$ is a translation-invariant invertible elliptic operator with the principal symbol $(\xi + ip)^{k_3}$;
\item[{\rm2)}] there exists a smooth homotopy of elliptic operators $A_\varepsilon$, $\varepsilon \in [0, 1]$, on the manifold $\mathcal{M}$ such that $A_0 = A \oplus 1_{N'}$, and the homotopy of the limit operators $\mathbf{a}_\varepsilon$ coincides with homotopy~\eqref{eq:homotopy2}.
\end{enumerate}
\end{proposition}
 
\begin{proof} 

1. The following symbol mapping for periodic $\psi$DOs:
$$
 \sigma_0:\Psi^0_{\mathrm{per}}(\mathbb{S}^1\times\mathbb{R})\longrightarrow  \mathrm{Mat}_N\bigl(C^\infty(S^*(\mathbb{S}^1_x\times \mathbb{S}^1_t) )\bigr)
$$ 
extends by continuity to the corresponding closures to $C^*$-algebras.
Furthermore, the following short exact sequence of $C^*$-algebras is defined:
\begin{equation}
\label{eq:short1}
  0\longrightarrow J\stackrel{j}\longrightarrow \overline{\Psi^0_{\mathrm{per}}(\mathbb{S}^1\times\mathbb{R})}\stackrel{\sigma_0}\longrightarrow  \mathrm{Mat}_N\bigl(C (S^*(\mathbb{S}^1_x\times \mathbb{S}^1_t) ) \bigr)\longrightarrow 0,
\end{equation}
where $J=\ker\sigma_0\subset \overline{\Psi^0_{\mathrm{per}}(\mathbb{S}^1\times\mathbb{R})}$ is the ideal of operators with zero symbol.
The Fourier transform on the real line yields the isomorphism $J \simeq C_0(\mathbb{R}, \mathcal{K}) \rtimes \mathbb{Z}$, where the group $\mathbb{Z}$ acts on the real line by integer shifts. In $K$-theory, sequence~\eqref{eq:short1} yields the following long exact sequence:
\begin{equation}
\label{eq:kth00}
  \ldots\longrightarrow K_1(J)\stackrel{j_*}\longrightarrow K_1(\overline{\Psi^0_{\mathrm{per}}(\mathbb{S}^1\times\mathbb{R})})\stackrel{\sigma_{0*}}\longrightarrow  K_1\bigl(C (S^*(\mathbb{S}^1_x\times \mathbb{S}^1_t) ) \bigr)
  \stackrel{\partial}\longrightarrow K_0(J)\longrightarrow \ldots
\end{equation}
Since the action of $\mathbb{Z}$ on $\mathbb{R}$ is proper and free, there exists the isomorphism
$$
 K_i(J)= K_i(C_0(\mathbb{R},\mathcal{K})\rtimes\mathbb{Z})\simeq K_i\bigl(C(\mathbb{R}/\mathbb{Z})\bigr)=K^i(\mathbb{S}^1)=\mathbb{Z},\qquad i=0,1.
$$
Under this isomorphism, the boundary map $\partial: K_1\bigl(C(S^*(\mathbb{S}^1_x \times \mathbb{S}^1_t))\bigr) \to \mathbb{Z}$ assigns to the class of a matrix-valued elliptic symbol on the torus the index of the corresponding elliptic operator on the torus.

2. Consider the symbol of the limit operator $\sigma_0(\mathbf{a})\in C^\infty\bigl(S^*(\mathbb{S}^1_x\times \mathbb{S}^1_t),{\rm GL}(N,\mathbb{C})\bigr)$ from the condition of the proposition being proved. We denote the class of this element in $K$-theory by
\begin{equation}
\label{eq:kth0}
[\sigma_0(\mathbf{a})]\in K^1\bigl(S^*(\mathbb{S}^1_x\times \mathbb{S}^1_t)\bigr).
\end{equation}
We have $S^*(\mathbb{S}^1_x\times \mathbb{S}^1_t)\simeq \mathbb{S}^1_x\times \mathbb{S}^1_t\times \mathbb{S}^1_{\xi,p}$. Hence, by the K\"unneth formula, we obtain
\begin{equation}
\begin{aligned}
\label{eq:kth1}
K^1\bigl(S^*(\mathbb{S}^1_x\times \mathbb{S}^1_t)\bigr) \simeq& K^1(\mathbb{S}^1_x\times \mathbb{S}^1_t)\otimes K^0(\mathbb{S}^1_{\xi,p})\oplus K^0(\mathbb{S}^1_x\times \mathbb{S}^1_t)\otimes K^1(\mathbb{S}^1_{\xi,p})  \\
\simeq& (\mathbb{Z}u\oplus \mathbb{Z}v)\otimes\mathbb{Z}{\bf 1}\oplus(\mathbb{Z}{\bf 1}\oplus\mathbb{Z}\beta)\otimes\mathbb{Z}w. 
\end{aligned}
\end{equation}
Here, the basis elements of the $K$-groups are as follows:
\begin{itemize}
\item $\bf 1$ is the class of one-dimensional trivial bundle;
\item $u=[e^{ix}]\in K^1(\mathbb{S}^1_x)$,  $v=[e^{it}]\in K^1(\mathbb{S}^1_t)$ , $w=[\xi+ip]\in K^1(\mathbb{S}^1_{\xi,p})$;
\item $\beta=[\gamma]-{\bf 1}\in K^0(\mathbb{S}^1_x\times \mathbb{S}^1_t) $  is the Bott generator (i.e., $\gamma$ is a one-dimensional bundle over a torus with the Chern number equal to one).
\end{itemize}
As follows from~\eqref{eq:kth1}, the element~\eqref{eq:kth0} can be expressed in terms of the constructed basis with integer coefficients $k_j$:
$$
[\sigma_0(\mathbf{a})]=k_1u\otimes {\bf 1}+k_2 v\otimes {\bf 1}+k_3 {\bf 1}\otimes w+ k_4 \beta\otimes w.
$$
Since the operator $A$ is assumed to be elliptic, its limit symbol $\mathbf{a}$ is invertible. In particular, the index of its Gelfand transform $\mathbf{a}(\theta)$ (see~\eqref{eq:atheta}) equals zero. On the other hand, the last expansion implies
\begin{equation}
\label{eq:ind0}
\ind\; \mathbf{a}(\theta)= k_1\ind (u\otimes {\bf 1})+k_2\ind( v\otimes {\bf 1})+k_3\ind( {\bf 1}\otimes w)+ k_4\ind( \beta\otimes w)= k_4.
\end{equation}
Here, the first three terms yield zero indices since the corresponding symbols do not depend on at least one of the variables. The index of the last term equals $1$, which can be easily verified using the Atiyah--Singer formula.
From~\eqref{eq:ind0}, we obtain $k_4=0$.
Using the explicit expression for a basis in the $K$-group, we arrive at the following equality in $K^1(S^*(\mathbb{S}^1_x\times \mathbb{S}^1_t))$:
\begin{equation*}
[\sigma_0(\mathbf{a})]=[e^{i(k_1x+k_2t)}\oplus(\xi+ip)^{k_3}].
\end{equation*}

3. Based on the symbol on the right-hand side in the last equality, we construct the following invertible elliptic $\psi$DO on a cylinder:
\begin{equation*}
\mathbf{a}_0=e^{i(k_1x+k_2t)}\oplus
    \left((-i\partial_x+\partial_t+1/2)(-\partial^2_x-\partial^2_t+1)^{-1/2}\right)^{k_3},
\end{equation*}
where the second term is a translation-invariant operator.
By construction, we obtain
\begin{equation}
\label{eq:smth1}
[\mathbf{a}]-[\mathbf{a}_0]\in \ker\sigma_{0*}={\rm im}\; j_*
\end{equation}
by the exactness of sequence~\eqref{eq:kth00}. It is easy to show that the embedding of the subalgebra
$$
J_0=C_0(\mathbb{R},\mathcal{K})\subset C_0(\mathbb{R},\mathcal{K})\rtimes\mathbb{Z}=J
$$
induces an isomorphism of $K_1$-groups. This and~\eqref{eq:smth1} imply that there exists an element $[\mathbf{a}'_0] \in K_1(C_0(\mathbb{R}, \mathcal{K}))$ for which the following equality holds:
\begin{equation*}
[\mathbf{a}]=[\mathbf{a}_0 \mathbf{a}'_0]\qquad\text{in group} \quad K_1(\overline{\Psi^0_{\mathrm{per}}(\mathbb{S}^1\times\mathbb{R})}).
\end{equation*}
This equality implies the existence of a stable homotopy between the invertible periodic elliptic operators $\mathbf{a}$ and $\mathbf{a}_0 \mathbf{a}'_0$. By smoothing this homotopy, we obtain the required homotopy~\eqref{eq:homotopy2}.

4. In part 2 of the proposition, the homotopy of elliptic operators on manifold $\mathcal{M}$ is constructed in the standard way.
Therefore, we omit this part of the proof.
\end{proof}

Let us return to the proof of Theorem~\ref{th:main'}.

1.  Let $A$ be an elliptic operator. By Proposition~\ref{prop1}, there exists a smooth stable homotopy from $A$ to $A_0$. Given the homotopy invariance of both the Fredholm index and the topological index, it suffices to verify their equality for the operator $A_0$.

2. Let us define an operator $A_{00}$ on $\mathcal{M}$ with the limit symbol (cf.~\eqref{eq:homotopy2})
$$
\mathbf{a}_{00}=e^{i k_1x}\oplus a (x,-i\partial_x,-i\partial_t)\oplus 1_{N'+N-2}
$$
and the interior symbol on $S^*\mathcal{M}$
$$
a_{00}(x,\xi,t,p)=\left\{
\begin{array}{ll}
 a_0(x,\xi,t,p) & \text{for }t\ge 0,\\
 e^{i k_1x+ik_2 t\chi(t)}\oplus (\xi+ip)^{k_3}\oplus 1_{N'+N-2} & \text{for }t\le 0.
\end{array}
\right.
$$
Here, the smooth cutoff function $\chi(t) \in C^\infty(\mathbb{R})$ is nondecreasing, with $\chi(t) \equiv 1$ for $t \ge 0$ and $\chi(t) \equiv 0$ for $t < -1/2$.

On the one hand, the following equality holds:
\begin{equation}
\label{eq:index51}
\ind A_0=\ind A_{00}.
\end{equation}
This equality follows from the fact that these two operators differ only on the set $\mathbb{S}^1\times (-\infty,0)$ and only in their first component, in which they are defined via bundle isomorphisms.
In particular, this modification does not alter the kernel or cokernel of the operator. The equality of the topological indices of these operators can be verified in an analogous manner:
\begin{equation}
\label{eq:index52}
\ind_t A_0=\ind_t A_{00}.
\end{equation}
On the other hand, the constructed operator $A_{00}$ has a translation-invariant limit symbol. Therefore, we can apply the index formula in~\cite{FST3} to this operator:
\begin{equation}
\label{eq:index53}
\ind  A_{00}=\ind_t A_{00}.
\end{equation}
In this case, the topological indices in the sense of the cited work and in the sense of the present work coincide.
 
Now, from~\eqref{eq:index51},~\eqref{eq:index52}, and~\eqref{eq:index53}, we obtain the desired equality between the analytic and topological indices of the operator $A_0$.

Theorem~\ref{th:main'} is proved.
Hence, Theorem~\ref{th:main} follows as a corollary.


\end{document}